\documentclass[12pt]{article}
\usepackage{bezier,latexsym,color}
\usepackage{epsfig,afterpage}
\usepackage{psfrag}
\usepackage{amsmath, amssymb}
\usepackage{bm}
\usepackage{color}
\usepackage{amsthm}
\usepackage{cite}
\usepackage[colorlinks = true,
            linkcolor = blue,
            urlcolor  = blue,
            citecolor = blue,
            anchorcolor = blue]{hyperref}
\newcommand\real{{\rm I\! R}}
\newcommand\nat{{\rm I\! N}}
\newcommand\cmplx{\;\hbox{\vrule height6.8pt width0.8pt depth-0.1pt
           \kern-3.6pt {\rm C}}}
\newcommand\beweisende{\ \hfill$\Box$\break}

\newtheorem{theorem}{Theorem}[section]
\newtheorem{definition}[theorem]{Definition}

\def\bql#1{\begin{eqnarray}\label{#1}}
\def\bq{\begin{equation}}
\def\eq{\end{equation}}
\def\eref#1{{\rm (\ref{#1})}}

\usepackage{array,multirow}
\usepackage{hhline}
{
\usepackage[table]{xcolor}

\begin{document}
\title{On the numerical solution of a hyperbolic initial boundary value problem by hypersingular boundary integral equations}

\author{Roman Chapko\thanks{Faculty of Applied Mathematics and
Informatics, Ivan Franko National University of Lviv, 79000
Lviv, Ukraine.} \, and Leonidas Mindrinos\thanks{Department of Natural Resources Development and Agricultural Engineering, Agricultural University of Athens, Greece.}}

 \date{}
 \maketitle

 \begin{abstract}
 In this study, we consider the numerical solution of the Neumann initial boundary value problem for the wave equation in 2D domains. Employing the Laguerre transform with respect to the temporal variable, we effectively transform this problem into a series of Neumann elliptic problems. The development of a fundamental sequence for these elliptic equations provides us with the means to introduce modified double layer potentials. Consequently, we are able to derive a sequence of boundary hypersingular integral equations as a result of this transformation. To discretize the system of equations, we apply the Maue transform and implement the Nystr\"om method with trigonometric quadrature techniques. To demonstrate the practical utility of our approach, we provide numerical examples.

\vspace*{1cm} {\em Keywords:} wave equation, Laguerre polynomials, double-layer potentials, hypersingular integral equations, Nystr\"om method.
\end{abstract}

\section{Introduction}

Hypersingular boundary integral equations frequently appear in different contexts, whether through the utilization of direct or indirect integral equation methodologies for addressing Neumann boundary value problems, as documented in \cite{WeHs, Vainikko1, Kr14}.  Notably, the representation of the solution as a double layer potential gives rise to a first-kind hypersingular boundary integral equation \cite{ChaJoSto, Kress2}. The motivation comes from the possibility to employ integral equations in the exploration and numerical solution of boundary value problems with Neumann boundary conditions on open boundaries \cite{ChaJoh, ChaKreMoe, Wendland}.

Various numerical methods have been developed to address such integral equations (see, for example, \cite{Vainikko1, Atkinson, KieKleeRath}). Typically, the Maue transform \cite{Maue, Mitzner} is employed prior to discretization, establishing a connection between the normal derivative of a double-layer potential and a single-layer potential. This yields a singular integro-differential equation and, in the context of two dimensions, results in a $2\pi$-periodic integro-differential equation, thereby making an effective utilization of the trigonometric Nystr\"om method possible \cite{Kress2}.

The principal objective of our paper is to employ hypersingular boundary integral equations in the context of time-dependent Neumann problems. This is accomplished in two steps. Firstly, we undertake temporal semi-discretization by employing the Laguerre transform, as detailed in \autoref{sec2a}. This process yields a sequence of 2D Neumann problems associated with elliptic equations. We derive the fundamental sequence governing these differential equations and introduce the corresponding modified double-layer potentials. This innovative approach enables the reduction of Neumann boundary value problems into a series of hypersingular boundary integral equations (BIE), see \autoref{sec2b}.

Subsequently, we apply the Maue transform to our modified potentials, resulting in a sequence of singular integro-differential equations of the first kind. The final phase involves a comprehensive discretization through the utilization of the Nystr\"om method, based on trigonometric quadratures, explained in \autoref{sec3}. In \autoref{sec4}, we present numerical examples to clarify our approach. It is worth noting that this two-step methodology has been successfully extended to address a diverse spectrum of initial boundary value problems for  parabolic and hyperbolic equations \cite{ChaJohMin20, ChaKre00, ChaMin19, ChaMin22}.

Prior to the implementation of our methodology, we provide the problem formulation. Let $\Omega$ be a bounded domain in $\real^2$ with a $C^3$ smooth boundary $\Gamma.$
 We consider the initial boundary value problem: 
\begin{equation}\label{wave1}
\frac{1}{c^2}\frac{\partial^2 u}{\partial t^2} =\Delta u, \quad \mbox{in } 
\Omega\times (0,\infty),
\end{equation}
together with homogeneous initial conditions
\begin{equation}\label{initial1}
\frac{\partial u}{\partial t}(\cdot\,,0)=u(\cdot\,,0)=0, \quad\mbox{in } \Omega,
\end{equation}
and the Neumann boundary condition
\begin{equation}\label{boundary1}
\frac{\partial u}{\partial \nu}= f, \quad \mbox{on }\Gamma\times[0,\infty).
\end{equation}
Here $c$ denotes the wave speed, $\nu$ is the outward unit normal vector to $\Gamma$, and $f$ is a given boundary function. 

\begin{theorem}
The initial boundary value problem  \eqref{wave1} -- \eqref{boundary1} is uniquely solvable.
\end{theorem}

\proof see \cite{Ladyzhenskaya1985}.
\beweisende

\section{The two-step approach}
\setcounter{equation}{0}
\setcounter{theorem}{0}

We solve numerically \eqref{wave1} -- \eqref{boundary1} in two steps. Initially we perform a semi-discretization for the time variable and then we apply the boundary integral equation method (BIE).

\subsection{Semi-discretization w.r.t. time variable}\label{sec2a}

We consider the following expansion
\begin{equation}\label{ut}
u(x,t) = \kappa \sum\limits_{n = 0}^\infty u_n(x) L_n(\kappa t),
\end{equation}
where
\begin{equation}
u_n (x) = \int_0^\infty e^{-\kappa t} L_n (\kappa t) u(x,t)dt, \quad n=0,1,2,\ldots,
\end{equation}
and $\kappa >0$ is a fixed parameter. Here, $L_n$ denote the Laguerre polynomials. Following \cite{ChaKre00, ChaMin22} we get the following sequence of stationary boundary value problems
\begin{equation}\label{eq_system1}
\Delta u_n - \gamma^2 u_n = \sum_{m=0}^{n-1} \beta_{n-m}u_m,  \quad   \mbox{in} \,\, D,
\end{equation}
for $n=0,1,\ldots,$ where $\beta_n = \tfrac{\kappa^2}{c^2} (n+1)$ and $\gamma^2=\beta_0.$ The boundary condition \eqref{boundary1} admits the form
\begin{equation}\label{eq_boundary2}
\frac{\partial u_n}{\partial \nu} = f_n,  \quad   \mbox{on} \,\, \Gamma, \end{equation}
where
\begin{equation}
f_n (x) = \int_0^\infty e^{-\kappa t} L_n (\kappa t) f(x,t)dt, \quad n=0,1,2,\ldots.
\end{equation}

\begin{theorem}\label{unique}
 The sequence of coupled problems \eqref{eq_system1} -- \eqref{eq_boundary2} has at most one solution.
 \end{theorem}
\proof 
By the Green's first theorem (see \cite{Kr14}) any solution
 $v\in C^2(D)\cap C(\bar D)$
 of $\Delta v- \gamma^2v=0$ in $D$ with homogeneous boundary
 condition $\frac{\partial v}{\partial \nu}=0$ on $\Gamma$ must vanish identically in $D$.
 Then the statement of the theorem follows by induction.
\beweisende

\subsection{The BIE method}\label{sec2b}

In order to establish the integral representation of the solution, it is necessary to first define the fundamental sequence. 

\begin{definition}\label{fund}
The sequence of functions $\left\{\Phi_n\right\}_{n=0}^{N}$ is called the fundamental sequence if it satisfies 
$$
\Delta_{x} \Phi_{n}(x, y)-\sum_{m=0}^{n} \beta_{n-m} \Phi_{m}(x, y)=\delta(x-y),
$$
where $\delta$ is the Dirac function.
\end{definition}

\begin{theorem}
 \label{theorem2}
Let $I_0$ and $I_1$ denote the modified Bessel functions \cite{AbSt}
\begin{equation}\label{besselmod}
 I_0(z)=\sum^\infty_{n=0} \;
 \frac{1}{(n!)^2}\,\left(\frac{z}{2}\right)^{2n},
 \quad
 I_1(z)=\sum^\infty_{n=0} \;
 \frac{1}{n!(n+1)!}\,\left(\frac{z}{2}\right)^{2n+1},
 \end{equation}
 and $K_0$ and $K_1$ the modified Hankel functions
 \begin{equation}\label{macdonald}
 \begin{aligned}
 K_0(z) &= - \left(\ln {z \over 2}+C\right)\,I_0(z)
 + \sum^\infty_{n=1}
 \frac{\psi(n)}{ (n!)^2} \,\left(\frac{z}{2}\right)^{2n},
 \\ 
 K_1(z) &= \frac{1}{z}+\left(\ln {z \over 2}+C\right)\,I_1(z)
 -\frac{1}{2}\sum^\infty_{n=0}
 \frac{\psi(n+1)+\psi(n)}{ n!(n+1)!}\,\left(\frac{z}{2}\right)^{2n+1}.
 \end{aligned}
 \end{equation} 
Here, we set $\psi(0)=0$,
 $$
 \psi(n)=\sum_{m=1}^n\frac{1}{m}\;, \quad n=1,2,\ldots,
 $$
 and let $C = 0.57721\ldots$
 denote Euler's constant.
 We define the polynomials
\begin{equation}\label{Polynomials}
 v_n (r)=
 \sum_{k=0}^{\left[\frac{n}{2}\right]}a_{n,2k} r^{2k},\quad
 w_n (r)=
 \sum_{k=0}^{\left[\frac{n-1}{2}\right]}a_{n,2k+1} r^{2k+1},
\end{equation}
 with the convention $w_0 (r) = 0,$ where 
\begin{align*}
a_{n,0} &=1, \\
 a_{n,n} &=-\frac{1}{2\gamma n}\;\beta_1 a_{n-1,n-1}, \\
a_{n,k} &=\frac{1}{2\gamma k}
 \left\{4\left[\frac{k+1}{2}\right]^2a_{n,k+1}
 -\sum_{m=k-1}^{n-1}\beta_{n-m} a_{m,k-1}\right\},
 \quad k=n-1,\ldots,1.
\end{align*}  
 
Then, the sequence of functions 
 \begin{equation}\label{fund-sol}
 \Phi_n (x, y) := \Phi_n (|x- y|)=K_0(\gamma |x-y|)v_n (|x-y|)+K_1(\gamma |x-y|)w_n (|x-y|),
 \end{equation}
 for $n=0,\ldots$, is a fundamental sequences of\/~\eref{eq_system1} in the sense of\/~{\rm Definition}~{\rm \ref{fund}}. 
\end{theorem}

\proof See  \cite{ChaKre00}. 
\beweisende

We represent the solution sequence of \eqref{eq_system1} -- \eqref{eq_boundary2} using a double-layer ansatz:
\begin{equation}\label{double}
 u_n(x)=\frac{1}{\pi} \sum_{m=0}^n
 \int_{\Gamma} \phi_m (y) \frac{\partial}{\partial \nu (y)}\Phi_{n-m}(x,y)\, ds(y),
 \quad x\in \Omega,
\end{equation}
with the unknown densities $\phi_m$ defined on the boundary $\Gamma$ and $\Phi_n$ is given by \eref{fund-sol}.

As revealed by the \autoref{theorem2} the fundamental sequence contains functions with logarithmic singularity, i.e.
\begin{equation}\label{PHI_expansion}
\Phi_n (r)\displaystyle =\ln \frac{1}{r}\left[I_0(\gamma r)v_n (r)-I_1(\gamma r)w_n (r)\right]+\ln\frac{2}{\gamma} - C + \frac{a_{n,1}}{\gamma}+O(r^2).
\end{equation}
Consequently, the modified double-layer potentials defined in \eref{double} exhibit the characteristics of a double layer potential for the Laplace equation. This leads us to the following result.
\begin{theorem}
The double-layer potential given by \eqref{double} solves the boundary value problems \eqref{eq_system1} -- \eqref{eq_boundary2} provided that the densities satisfy the following sequence of boundary integral equations
\begin{equation}\label{system1}
(\mathcal{T}_0 q_n ) (x) = f_n (x)-\sum\limits_{m=0}^{n-1}(\mathcal{T}_{n-m} q_m ) (x), \quad   x\in\Gamma, \quad n = 0,1,\ldots,
\end{equation}
where $\mathcal{T}_n$ is the normal derivative of the double layer operator, namely
\begin{equation}\label{normal_double_layer}
(\mathcal{T}_n f ) (x) = \frac{1}{\pi} \frac{\partial}{\partial \nu (x)} \int_{\Gamma}f(y) \frac{\partial }{\partial \nu (y)}  \Phi_n (x,y)  ds(y),\quad x\in\Gamma.
\end{equation}
\end{theorem}

In the following theorem we derive a Maue type formula that reduces the hypersingularity of the operator $\mathcal{T}_n,$ see for example the classic works \cite{Mitzner, Maue}.

\begin{theorem}\label{Maue}
The following decomposition holds true:
\begin{equation}\label{maue}
\frac{\partial}{\partial \nu (x)} \sum_{m=0}^n(\mathcal{D}_{n-m} q_m ) (x) = \sum_{m=0}^n \left[(\mathcal{H}_{n-m} q_m ) (x)-(\mathcal{L}_{n-m} q_m ) (x)\right], \quad n = 0,1,\ldots
\end{equation}
where
\[
\begin{aligned}
(\mathcal{H}_{n} f ) (x) &=\frac{1}{\pi}\frac{d}{d s(x)}\int_{\Gamma}  \frac{d f}{d s}(y) \Phi_n (x,y)ds(y),\\
(\mathcal{L}_{n} f ) (x) &=\frac{1}{\pi}\nu(x)\cdot \int_{\Gamma}f(y) \nu(y)\sum_{k=0}^n\beta_{n-k}\Phi_k(x,y)ds(y).
\end{aligned}
\]
\end{theorem}

\proof The formula \eqref{maue} is the result of lengthy but straightforward calculations, for $x\in\Gamma.$
\beweisende



The evaluation of the operator $\mathcal{H}_n$ requires the derivatives of $\Phi_n$ which can be obtained using the recurrent relations of the modified Hankel functions. We obtain
\begin{equation*}
\begin{aligned}
\Psi_n(r) &:=\Phi_n^{\prime}(r)=K_0(\gamma r)\tilde{v}_n (r)+K_1(\gamma r)\tilde{w}_n (r), \\
\Omega_n (r) & :=\Phi_n^{\prime \prime}(r)=K_0(\gamma r){\bar v}_n (r)+\frac{K_1(\gamma r)}{r}{\bar w}_n (r), 
\end{aligned}
\end{equation*}
for the polynomials
\begin{align*}
\tilde{v}_n(r) &=2\sum_{k=1}^{\left[\frac{n}{2}\right]}ka_{n,2k}r^{2k-1}-\gamma w_n (r), \\
\tilde{w}_n (r) &=2\sum_{k=1}^{\left[\frac{n-1}{2}\right]}ka_{n,2k+1} r^{2k}-\gamma v_n (r), \\
\bar{v}_n (r) &=2\sum_{k=1}^{\left[\frac{n}{2}\right]}k(2k-1)a_{n,2k} r^{2k-2}-\gamma \sum_{k=0}^{\left[\frac{n-1}{2}\right]}(2k+1)a_{n,2k+1} r^{2k}-\gamma \tilde{w}_n(r), \\
\bar{w}_n (r) &=4\sum_{k=1}^{\left[\frac{n-1}{2}\right]}k^2a_{n,2k+1} r^{2k}-2\gamma \sum_{k=1}^{\left[\frac{n}{2}\right]}ka_{n,2k} r^{2k}-\gamma r\tilde{v}_n (r)- \tilde{w}_n (r),
\end{align*}
with $\bar{v}_0 (r)=\gamma^2.$ Using the series expansion of the modified Bessel functions, it can be inferred that
\begin{align}
\label{PSI_expansion}
\Psi_n (r)&\displaystyle =-\frac{1}{r}+\ln \frac{1}{r}\left[I_0(\gamma r)\tilde{v}_n (r)-I_1(\gamma r)\tilde{w}_n (r)\right] +O(r),\\
\Omega_n (r)&\displaystyle =\frac{1}{r^2}+\ln \frac{1}{r}\left[I_0(\gamma r)\bar{v}_n (r)-\frac{I_1(\gamma r)}{r}\bar{w}_n (r)\right]+\omega_n +O(r^2),
\end{align}
with
$$
\omega_n =\left(\ln\frac{2}{\gamma} - C\right)\left(\frac{\gamma^2}{2}-\gamma  a_{n,1}+2a_{n,2}\right) -\frac{\gamma^2}{4}+\gamma a_{n,1} - 3a_{n,2}+\frac{2a_{n,3}}{\gamma}.
$$

\section{The Nystr\"om method}\label{sec3}
\setcounter{equation}{0}
\setcounter{theorem}{0}

In this section, we address the numerical solution of \eqref{system1} through the application of the Nystr\"om method. Singular kernels are managed using quadrature rules.

\subsection{Boundary parametrization}
Let the boundary curve $\Gamma\in C^3,$ have the parametric representation 
$$
\Gamma=\left\{x(s)=\left(x_{1}(s), x_{2}(s)\right), s \in[0,2 \pi]\right\},
$$
where $x: \mathbb{R} \rightarrow \mathbb{R}^{2}$ is a $2 \pi$-periodic with $\left|x^{\prime}(s)\right|>0,$ for all $s$. We define $r(s,\sigma)=|x(s)-x(\sigma)|.$ 

We consider the parametrized forms of the integral operators for $n=0,1,\ldots$
\begin{align*}
(H_n f ) (s) & =\frac{1}{2\pi}  \int_0^{2\pi} f^\prime (\sigma)Q_{n}(s,\sigma) d\sigma, \\
(L_n f ) (s) & =\frac{1}{2\pi}  \int_0^{2\pi} f (\sigma) \tilde{W}_{n} (s,\sigma) d\sigma,
\end{align*}
for the kernels
\begin{align*}
W_n (s,\sigma) &=2 |x^\prime (\sigma)|\Phi_{n} (x (s),x (\sigma)), \\
 \tilde{W}_{n}(s,\sigma) &=\nu(x(s))\cdot \nu(x(\sigma)) \sum_{k=0}^n \beta_{n-k} W_{k}(s,\sigma), \\
 Q_{n}(s,\sigma) &=\frac{2}{|x^\prime(s)|}\frac{\partial }{\partial s}\Phi_{n} (x (s),x (\sigma)) = 2\Psi_{n} \left(r(s,\sigma)\right) h (s,\sigma),
\end{align*}
for $s \neq \sigma,$ and $n=1, \ldots, N$. 
Here we have introduced the function
\[
h (s,\sigma) = \frac{(x(s)-x(\sigma))\cdot x^\prime (s)}{|x^\prime(s)|r(s,\sigma)}.
\]

Then, the integral equation \eqref{system1} admits the parametric form
\begin{equation}\label{system2}
((H_0 - L_0) \psi_n ) (s) = f_n (x(s)) - \sum\limits_{m=0}^{n-1}((H_{n-m} - L_{n-m})\psi_m ) (s),
\end{equation}
for $s\in[0,2\pi]$ and $n=0, \ldots, $ with $\psi_{n}(s)=q_{n}\left(x(s)\right).$

\subsection{Singularities}
Given the expansions \eref{PHI_expansion} and \eref{PSI_expansion}, it becomes apparent that the kernels exhibit logarithmic singularity.
The kernel of the single layer operator admits the decomposition
$$
W_{n}(s,\sigma)=W^1_{n}(s,\sigma)\ln\left(\frac{4}{e}\sin\frac{s-\sigma}{2}\right) +W^2_{n}(s,\sigma),
$$
where
$$
W^1_{n}(s,\sigma)=-I_0(\gamma r(s,\sigma))v_n(r(s,\sigma))+I_1(\gamma r (s,\sigma))w_n(r(s,\sigma))]
$$ 
and
$$
W^2_{n}(s,\sigma)=\left\{
\begin{array}{cc}
\displaystyle W_{n}(s,\sigma)-W^1_{n}(s,\sigma)\ln\left(\frac{4}{e}\sin\frac{s-\sigma}{2}\right) & \mbox{for}\; s\ne\sigma,\\
\displaystyle \left[-2\ln\left(\frac{\gamma |x^\prime(s)|}{2}\right)-1-2C+\frac{2a_{n,1}}{\gamma}\right] |x^\prime(s)|  & \mbox{for}\; s=\sigma.
\end{array}
\right.
$$

To handle the singularity of the operator $L_n,$ we apply the following transformation:
$$
\frac{1}{2\pi}  \int_0^{2\pi} f^\prime (\sigma)Q_{n}(s,\sigma) d\sigma 
 = \frac{1}{2\pi|x^\prime(s)|}  \int_0^{2\pi} f^\prime (\sigma)\cot \frac{\sigma-s}{2}d\sigma +\frac{1}{2\pi}  \int_0^{2\pi} f (\sigma)\tilde{Q}_{n}(s,\sigma) d\sigma,
$$
where
$$
\tilde{Q}_{n}(s,\sigma) 
= 2\Omega_n (r(s,\sigma)) h_1(s,\sigma) + 2\Psi_n (r (s,\sigma)) h_2(s,\sigma) + \frac{1}{2|x^\prime(s)|\sin^2\frac{s-\sigma}{2}},
$$
for the functions
\begin{align*}
h_1(s,\sigma) &=\frac{(x(\sigma)-x(s))\cdot x^\prime(\sigma)(x(s)-x(\sigma))\cdot x^\prime(s)}{|x^\prime(s)|r^2(s,\sigma)}, \\
h_2(s,\sigma) &=-\frac{x^\prime(s)\cdot x^\prime (\sigma)}{|x^\prime(s)|r (s,\sigma)}+\frac{(x (s)-x (\sigma))\cdot x^\prime(\sigma)(x(s)-x(\sigma))\cdot x^\prime(s)}{|x^\prime(s)| r^3(s,\sigma)}.
\end{align*}
Note that $\lim_{\sigma\to s} h_1(s,\sigma)=-|x^\prime(s)|,$ and $\lim_{\sigma\to s} h_2(s,\sigma)=0$. Again, a straightforward calculation yields the following representation
$$
\tilde{Q}_{n}(s,\sigma) =\tilde{Q}_{n}^1(s,\sigma)\ln\left(\frac{4}{e}\sin\frac{s-\sigma}{2}\right) + \tilde{Q}_{n}^2(s,\sigma) 
$$
with
\begin{align*}
 \tilde{Q}_{n}^1(s,\sigma) &=-h_1(s,\sigma)\left[I_0(\gamma r(s,\sigma))\bar{v}_n (r(s,\sigma))-\frac{I_1(\gamma r(s,\sigma))}{r}\bar{w}_n (r(s,\sigma))\right]\\
&\phantom{=}-h_2(s,\sigma)\left[I_0(\gamma r (s,\sigma))\tilde{v}_n (r(s,\sigma))-I_1(\gamma r(s,\sigma))\tilde{w}_n(r(s,\sigma))\right]
\end{align*}
and
$$
\tilde{Q}_{n}^2(s,\sigma)=\left\{
\begin{array}{cc}
\displaystyle \tilde{Q}_{n}(s,\sigma) -\tilde{Q}_{n}^1(s,\sigma)\ln\frac{4}{e}\sin\frac{s-\sigma}{2} & \mbox{for}\; s\ne\sigma,\\
\displaystyle -\ln(e|x^\prime(s)|^2)\tilde{Q}_{n}^1(s,s)+2 \omega_n h_1(s,s)+\frac{1}{6|x^\prime(s)|} \\
\displaystyle +\frac{(x^\prime(s)\cdot x^{\prime\prime}(s))^2}{|x^\prime(s)|^5} -\frac{x^\prime(s)\cdot x^{(3)}(s)}{3|x^\prime(s)|^3}-\frac{x^{\prime\prime}(s)^2}{2|x^\prime(s)|^3} & \mbox{for}\; s=\sigma.
\end{array}
\right.
$$
Recall that $\tilde{Q}_{n}^1(s,s)=-(\frac{\gamma^2}{2}-\gamma a_{n,1}+2a_{n,2}) h_1(s,s)$.


\begin{theorem}\label{well-possed}
For any $f_n\in C^{0,\alpha}[0,2\pi],$ the sequence of equations \eqref{system2} admits unique solutions $\psi_n\in C^{1,\alpha}[0,2\pi]$.
\end{theorem}
\proof 
For $n=0,$ we can rewrite the first equation of \eqref{system2} in the following equivalent form
$$
(\tilde{H}_0 - \tilde{L}_0)\psi_0 (s) =f_0(s),
$$
where 
$$
\tilde{H}_0\phi(s)=\frac{1}{2\pi}\int_0^{2\pi}\left(\phi^\prime(\sigma) \cot\frac{\sigma-s}{2} + \phi(\sigma)\right)d\sigma,
$$
and 
$$
\tilde{L}_0\phi(s)=\frac{1}{2\pi}\int_0^{2\pi}\phi(\sigma)\left(\tilde{Q}_0(s,\sigma)-1\right)d\sigma.
$$
The operator $\tilde{H}_0: C^{1,\alpha}[0,2\pi]\to C^{0,\alpha}[0,2\pi]$ is bounded and has a bounded inverse and the operator $\tilde{L}_0: C^{1,\alpha}[0,2\pi]\to C^{0,\alpha}[0,2\pi]$ is compact (see for more details \cite{Kress2}). Thus, from the Riesz theory
for compact operators we have that $\tilde{H}_0 - \tilde{L}_0: C^{1,\alpha}[0,2\pi]\to C^{0,\alpha}[0,2\pi]$  has a bounded inverse if
and only if $\tilde{H}_0 - \tilde{L}_0$ is injective.  The injectivity of this operator follows from the uniqueness \autoref{unique}. The statement of the theorem follows by induction.
\beweisende

\subsection{Full discretization}
We utilize the following standard quadrature rules for the estimation of the integrals associated with the equation \eqref{system2} \cite{Kr14}
\begin{align*}
 \frac{1}{2\pi}
 \int_0^{2\pi}
 f(\sigma)\,
 d\sigma
 &\approx \frac{1}{2M}
 \sum_{k=0}^{2M-1}
  f(s_k), \\
   \frac{1}{2\pi}
 \int_0^{2\pi}
 f(\sigma)
 \ln\left(\frac{4}{e} \sin^2 {s- \sigma \over 2}\right)
 d\sigma
 &\approx
 \sum_{k=0}^{2M-1}
 {R}_{k}(s)\,f(s_k), \\
  \frac{1}{2\pi}
 \int_0^{2\pi}
 f^\prime (\sigma)
 \cot {\sigma -s  \over 2}
 d\sigma
 &\approx
 \sum_{k=0}^{2M-1}
 {T}_{k}(s)\,f(s_k),
\end{align*}
for the mesh points 
\begin{equation}\label{mesh_points}
s_k=kh, \quad k=0,\ldots,2M-1, \quad h=\pi/M,
\end{equation} 
and the weight functions 
\begin{align*}
 {R}_{k}(s) &= 
\displaystyle
 - {1 \over 2M} \;  \left(1+2\sum^{M-1}_{m=1} \; {1 \over m} \,
\cos m(s- s_k) - {1 \over n} \, \cos M(s-s_k )\right), \\
 {T}_{k}(s) &= 
\displaystyle
 - {1\over M} \;  \sum^{M-1}_{m=1} \; m
\cos m(s- s_k) - {1\over 2} \, \cos M(s-s_k ).
\end{align*}

Subsequently, we proceed to impose collocation upon the resultant approximate equations, leading to the derivation of a sequence of linear systems
\begin{equation}\label{system3}
\sum_{j=0}^{2M-1} M_{0;ij}\psi_{n;j}  = f_n(x(s_i)) -\sum_{m=0}^{n-1}\sum_{j=0}^{2M-1} M_{n-m;ij} \psi_{m;j},
\end{equation}
for $i=0,\ldots,2M-1,$ and the matrix-valued coefficient
\[
M_{n;ij} =\frac{T_j(s_i)}{|x^\prime(s_i)|}+\left[\tilde{Q}_n^1(s_i,s_j)-\tilde{W}_n^1(s_i,s_j)\right]R_j(s_i)+\frac{1}{2M}\left[\tilde{Q}_n^2(s_i,s_j)-\tilde{W}_n^2(s_i,s_j)\right].
\]

For the numerical solution $\tilde{\psi}_n$, which is constructed by trigonometrical interpolation using the values $\psi_{n;j}$ from \eqref{system3}  we have the following error estimate (see for details \cite{Kress2})
$$
\|\tilde{\psi}_n - \psi_n \|_{1,\alpha}\le c_n\frac{\ln M}{M^{q+\beta-\alpha}}\|f_n\|_{q,\beta}
$$
for $0<\alpha\le\beta<1$ and $q\in \nat$.

Then, the numerical solution of the sequence of stationary problems \eqref{eq_system1} -- \eqref{eq_boundary2} following the representation \eqref{double} admits the following form
\begin{equation}\label{sol_stationary}
u_n (x)\approx u_{n,M} (x)=\frac{1}{M}\sum_{m=0}^n\sum_{j=0}^{2M-1} \psi_{m,j}\Psi_{n-m}(x,x (s_j)) h (x,s_j), \quad x\in D.
\end{equation}

In conclusion, the numerical approximation of the initial boundary problem, as defined in \eqref{wave1} -- \eqref{boundary1}, employing the representation \eqref{ut}, yields the following approximate solution:
\begin{equation}\label{sol_dynamic}
u(x,t)\approx u_{N,M}(x,t)=\kappa\sum_{n=0}^N u_{n,M}(x)L_n(\kappa t), \quad (x,t)\in D\times (0,T]. 
\end{equation}

\section{Numerical results}\label{sec4}
\setcounter{equation}{0}
\setcounter{theorem}{0}

We illustrate the effectiveness of our approach through the examination of different numerical examples. Initially, we aim to provide numerical validation for the formulated equation \eqref{system3}. To achieve this, we choose in the first example an arbitrary source point denoted as $z_1 \in \real^2 \setminus D$ and proceed to construct the boundary function:
\[
f_n (x(s)) = \Psi_n (|x(s) - z_1|) \frac{(x(s)- z_1)\cdot \nu (x(s))}{|x(s)- z_1|}.
\]
We know the exact solution of the stationary problem at a given interior point $z_2 \in D,$ given by $u^{ex}_n (z_2) = \Phi_n (|z_2 - z_1|),$ and we compare it with the numerical one, see \eqref{sol_stationary}. We expect to achieve exponential convergence with respect to the number of grid boundary points $M$ \cite{Kr14}.

\begin{figure}[t]
\begin{center}
\includegraphics[scale=0.9]{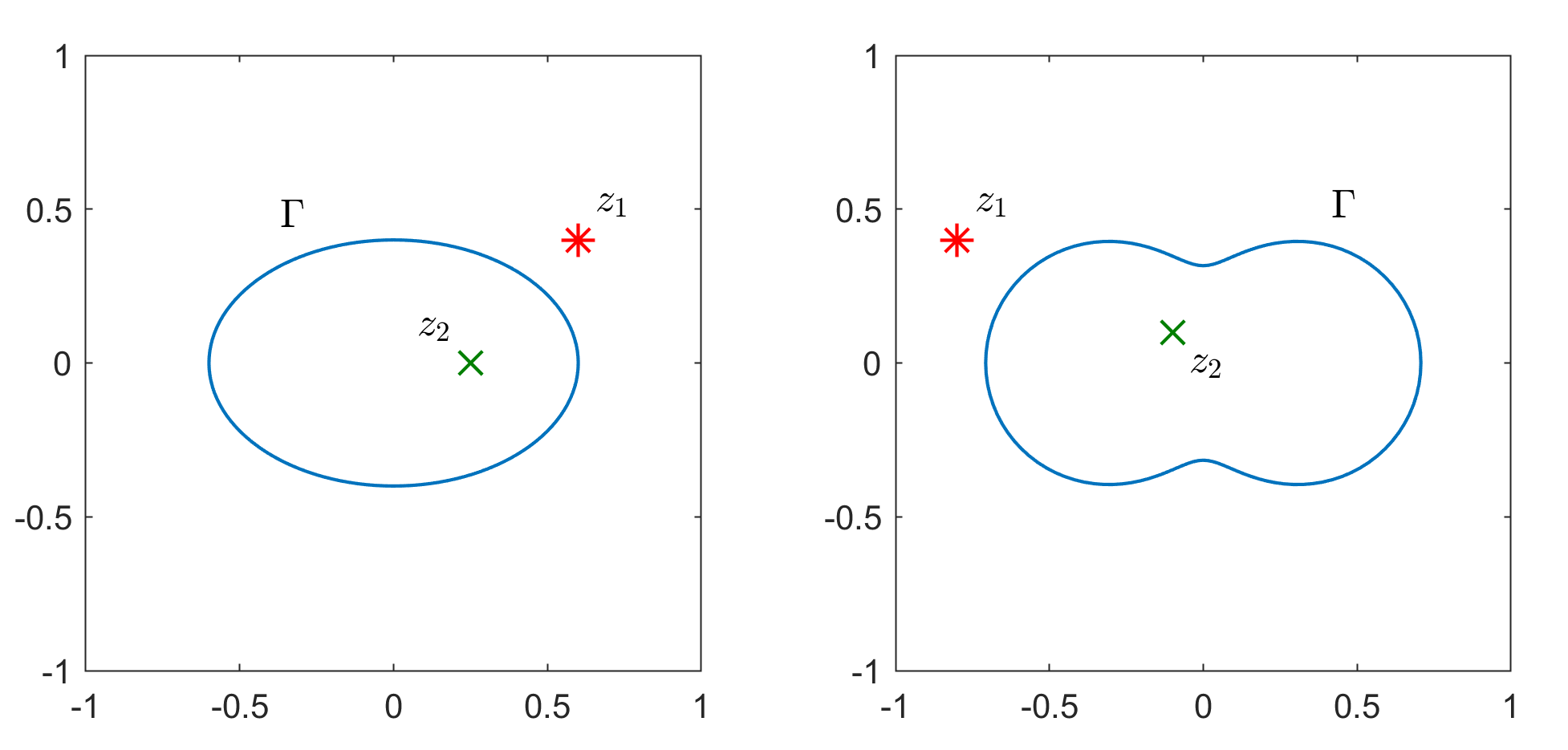}
\caption{The geometry of the problem. The boundary curve (solid blue line), the source point (red star) and the measurement point (green cross) for the first (left) and the second (right) example. }\label{fig1}
\end{center}
\end{figure}

The parametrized boundary curve is given by
\[
x(s) = (0.6 \cos (s), \, 0.4 \sin (s)), \quad s\in [0,2\pi],
\]
and we set $\kappa = c = 1.$ In Table \ref{table1} we present the two solutions for $n=0,1,5,$ at the measurement point $z_2 = (0.25, 0),$ for a source point at $z_1 = (0.6, 0.4),$ see the left picture in \autoref{fig1}.

\setlength{\tabcolsep}{18pt} 
\begin{table}[t]
\begin{center}
 \begin{tabular}{| c  | c  | c  | c  | } 
 \hline
 $M$ & $ u_0 (z_2) $ & $u_1 (z_2) $ & $u_5 (z_2) $  
\\ \hline 
8 &    0.8719916644   &  0.0616968878  &  $ -0.3179550968 $    \\ 
16 &  0.8740840870  &   0.0607702205 &   $  -0.3184694131$  \\
32 &    0.8742487939    & 0.0606328953   &  $-0.3184787687 $  \\
64 &   0.8742487910  &   0.0606328977  &   $-0.3184787686$   \\
 \cline{1-1}\hhline{~===}

 \multicolumn{1}{c|}{} &  $ u_0^{ex} (z_2) $   & $ u_1^{ex} (z_2) $ & $ u_5^{ex}  (z_2 ) $ \\
 \cline{2-4} 
\multicolumn{1}{c|}{} &       0.8742487910  &   0.0606328977    & $-0.3184787686 $
 \\ \cline{2-4}
\end{tabular}
\vspace{0.4cm}
\caption{The computed and the exact solutions of \eref{eq_system1}--\eref{eq_boundary2} of the first example.}\label{table1}
\end{center}
\end{table}

In the second example, we test the applicability of the proposed scheme  for different boundary curve and parameters. The geometry of the problem is presented in the right picture of \autoref{fig1}. We consider a peanut-shaped boundary of the form
\[
x(s) = \sqrt{ 0.5 \cos^2 (s) + 0.1 \sin^2 (s)}
(\cos (s), \,  \sin (s)), \quad s\in [0,2\pi],
\]
for $\kappa = 1$ and $c=1/2.$ The source point is now located in $z_1 = (-0.8,\, 0.4),$ and we measure at $z_2 = (-0.1,\, 0.1).$ The computed values are presented in Table \ref{table2}.

\setlength{\tabcolsep}{18pt} 
\begin{table}[t]
\begin{center}
 \begin{tabular}{| c  | c  | c  | c  | } 
 \hline
 $M$ & $ u_0 (z_2) $ & $u_1 (z_2) $ & $u_5 (z_2) $  
\\ \hline 
8 &    0.2313220521    & $-0.2317638260$  &   0.0887619852     \\ 
16 &  0.2078149403 &    $-0.2016053322$ &    0.0867268360  \\
32 &    0.2074883363    & $-0.2012217392$ &    0.0867114468   \\
64 &    0.2074882774 &     $-0.2012216749$   &  0.0867114477   \\
 \cline{1-1}\hhline{~===}

 \multicolumn{1}{c|}{} &  $ u_0^{ex} (z_2) $   & $ u_1^{ex} (z_2) $ & $ u_5^{ex}  (z_2 ) $ \\
 \cline{2-4} 
\multicolumn{1}{c|}{} &       0.2074882774  &   $-0.2012216749$  &    0.0867114477 
 \\ \cline{2-4}
\end{tabular}
\vspace{0.4cm}
\caption{The computed and the exact solutions of \eref{eq_system1}--\eref{eq_boundary2} of the second example.}\label{table2}
\end{center}
\end{table}

\begin{figure}[t]
\begin{center}
\includegraphics[scale=0.9]{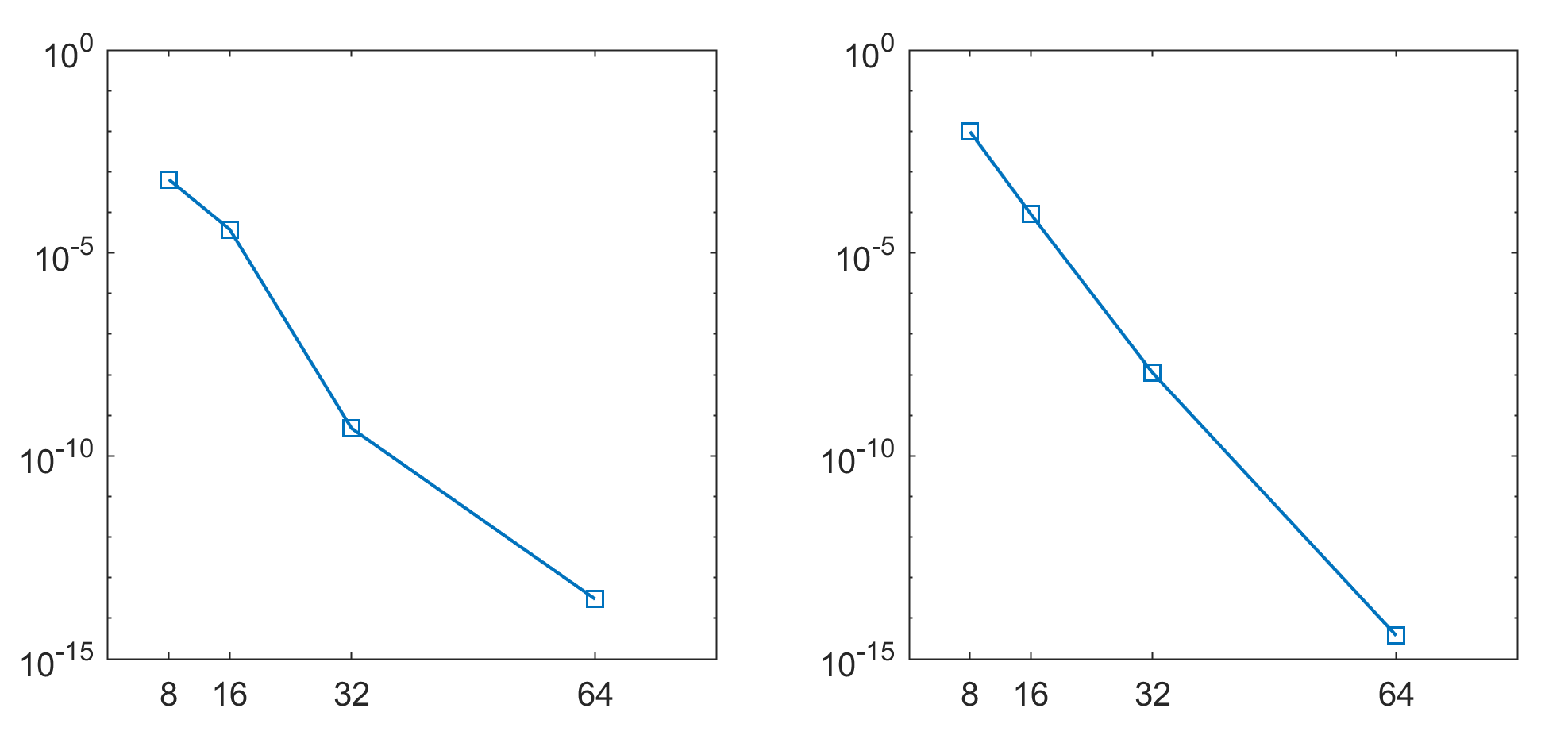}
\caption{The $L^2$ norm of the difference between the computed and the exact solutions, in semi-logarithmic scale, of the first (left) and the second (right) example. }\label{fig2}
\end{center}
\end{figure}

Exponential convergence with respect to the spatial discretization $M$ is readily evident, as can be observed in \autoref{fig2}, which depicts a semi-logarithmic representation of the $L^2$ norm of the difference between the exact and numerical solution. 

In the third example we examine the solvability of the problem \eqref{wave1} -- \eqref{boundary1} at different time steps. We consider the boundary function
\[
f(t) =  \frac{t^2}{4}e^{-t+2},
\]
which admits the expansion
\[
f(t) = \frac{\kappa e}{4} \sum_{n=0}^\infty \frac{2+\kappa n (\kappa (n-1)-4)}{(\kappa+1)^{n+3}} L_n (\kappa t).
\]

We set $\kappa=1,$ and $c=1/2,$ and the boundary is kite-shaped with parametrization
\[
x(s) = (0.2\cos (s) + 0.1 \cos (2s) - 0.2, \,  0.2\sin (s)+0.1), \quad s\in [0,2\pi],
\]

In this particular case, the absence of an exact solution prevent us from a direct comparison; however, we notice the convergence concerning both the parameter $M$ and the count of Fourier coefficients $N$ through the data presented in Table \ref{table3}.

\setlength{\tabcolsep}{18pt} 
\begin{table}[ht]
\begin{center}
 \begin{tabular}{| c | c  | c  | c  | c  | } 
 \hline
$t$ & $M$ & $N=10$ & $N=15 $ & $N=20 $  
\\ \hline 
 \multirow{3}{*}{1} & 16 & 0.0834968169 & 0.0857192560 & 0.0858032770 
 \\
 & 32 &  0.0835172830 & 0.0857399022 & 0.0858239138
 \\ 
  & 64 &  0.0835172976 & 0.0835172976 & 0.0858239286
  \\ 
 \cline{1-1}\hhline{=====}
  \multirow{3}{*}{2} & 16 & 0.8426917166  & 0.8398099958 &  0.8396857522
 \\
 & 32 & 0.8428312210 &  0.8399492790  & 0.8398249923
  \\ 
  & 64 & 0.8428313185 & 0.8399493762   & 0.8398250895
  \\ 
 \cline{1-1}\hhline{=====}
  \multirow{3}{*}{3} & 16 &  2.6020748549 & 2.6040387055 &   2.6043197336
 \\
 & 32 &   2.6023679039  &  2.6043317618  & 2.6046127580
 \\ 
  & 64 & 2.6023681749  &  2.6043320330 &  2.6046130293 
 \\ \hline
\end{tabular}
\caption{The numerical solution at $z_2 = (-0.2,\, 0.2)$ of the problem \eqref{wave1}--\eqref{boundary1} for the setup of the third example.}\label{table3}
\end{center}
\end{table}

\begin{figure}[t]
\begin{center}
\includegraphics[scale=0.9]{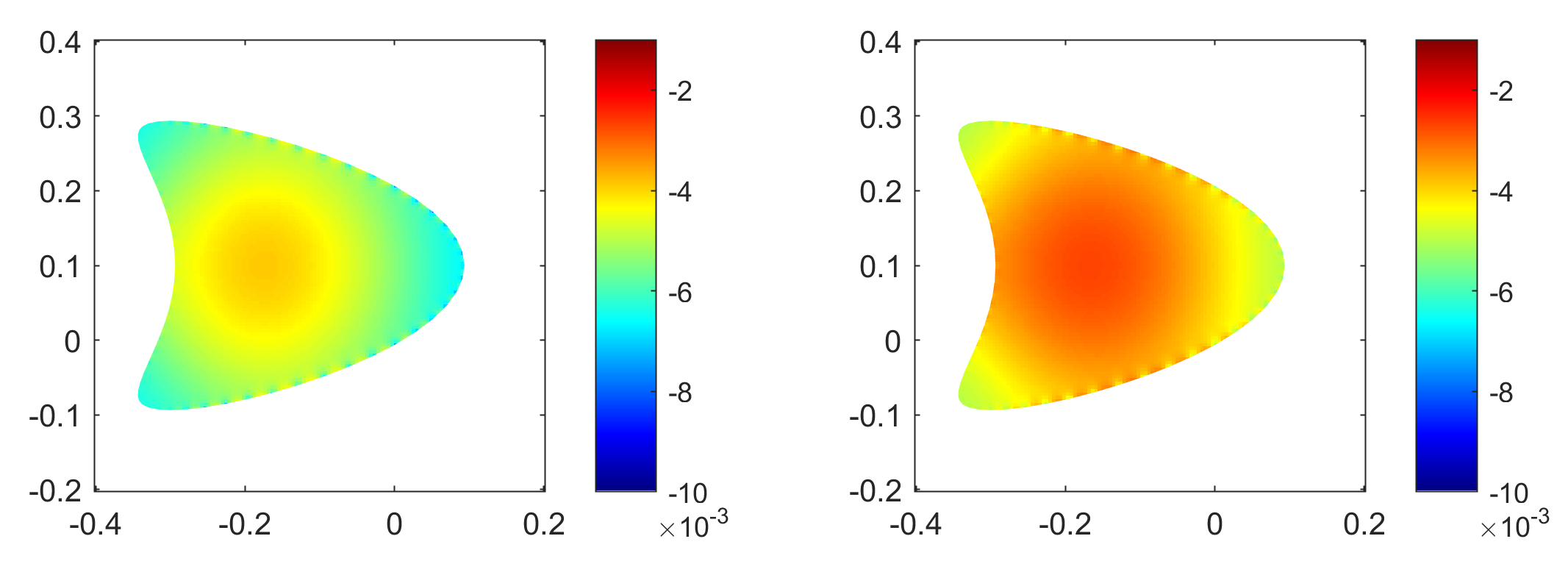}
\caption{The computed solution $u_{10,32},$ given by \eqref{sol_dynamic}, for $t=0.001$ (left) and $t=0.003$ (right) for the setup of the third example. }\label{fig3}
\end{center}
\end{figure}

In \autoref{fig3}, we present the numerical solution of the third example for $N=10$ and $M=32,$ at the time steps $t=0.001$ and $t=0.003.$

\section{Conclusions}
In this work, we have introduced a hypersingular integral equation approach for hyperbolic initial boundary value problems with Neumann boundary conditions. Our methodology is centered around a systematic dimensionality reduction strategy. By applying the Laguerre transform, we reformulate the time-dependent problem into a sequence of stationary Neumann problems for inhomogeneous elliptic equations. Subsequently, we construct modified double layer potentials, enabling us to derive a sequence of hypersingular boundary integral equations. Consequently, the originally three-dimensional problem undergoes a transformation, reducing to a sequence of one-dimensional integral equations.  The full discretization is achieved through the application of the trigonometric quadrature method, making this methodology suitable for both interior and exterior domains. Our numerical experiments serve as a proof of the practicality and efficacy of the methodology we have developed.

This work fills a gap in the field of hyperbolic boundary value problems, offering a versatile and adaptive tool for tackling such challenges. The dimension reduction, an  important feature of our approach, not only simplifies the computational complexity but also preserves accuracy. Future research can delve into expanding this approach to a broader class of problems and more complex geometries, including scenarios involving open boundaries.

\bibliographystyle{siam}
\bibliography{chamin_refs}

 \end{document}